\documentclass[11pt,fleqn,twoside]{article}
\usepackage{amsfonts,amssymb,latexsym}
\makeatletter
\newcommand{\prava}[1]{\small\it
\begin{flushleft}
Copyright \copyright \ 1999 by  #1
\end{flushleft}}

\newcommand{\name}[1]{\begin{flushleft}
                       \LARGE \bf #1
                       \end{flushleft}\vspace{-3mm}}

\newcommand{\Author}[1]{\begin{flushleft}
                       \it #1 \end{flushleft}}

\newcommand{\Adress}[1]{\begin{flushleft}
                       \it #1 \end{flushleft}}

\newcommand{\Date}[1]{\begin{flushleft}
                      \small  \it #1 \end{flushleft}}

\newcommand{\ehkol}{Author \ name}
\newcommand{\ohkol}{Article \ name}
\renewcommand{\@evenhead}{
\hspace*{-3pt}\raisebox{-15pt}[\headheight][0pt]{\vbox{\hbox to \textwidth 
{\thepage \hfil \ehkol}\vskip4pt \hrule}}}
\renewcommand{\@oddhead}{
\hspace*{-3pt}\raisebox{-15pt}[\headheight][0pt]{\vbox{\hbox to \textwidth 
{\ohkol \hfil \thepage}\vskip4pt\hrule}}}
\renewcommand{\@evenfoot}{}
\renewcommand{\@oddfoot}{}

\newcommand{\be}{\begin{equation}}
\newcommand{\ee}{\end{equation}}
\newcommand{\ba}{\hspace*{-5pt}\begin{array}}
\newcommand{\ea}{\end{array}}

\newcommand{\ds}{\displaystyle}
\makeatother

\newtheorem{theorem}{Theorem}[section]
\newtheorem{remark}{Remark}[section]
\newtheorem{proposition}{Proposition}[section]

\begin{document}

\thispagestyle{empty}
\setcounter{page}{344}
\renewcommand{\ehkol}{B.A. Kupershmidt and O.S. Stoyanov}
\renewcommand{\ohkol}{Coadjoint Poisson Actions of Poisson-Lie Groups}

\begin{flushleft}
\footnotesize \sf
Journal of Nonlinear Mathematical Physics \qquad 1999, V.6, N~3,
\pageref{kupershmidt_5-fp}--\pageref{kupershmidt_5-lp}.
\hfill {\sc Article}
\end{flushleft}

\vspace{-5mm}

\renewcommand{\footnoterule}{}
{\renewcommand{\thefootnote}{}
 \footnote{\prava{B.A. Kupershmidt and O.S. Stoyanov}}}

\name{Coadjoint Poisson Actions of Poisson-Lie Groups}\label{kupershmidt_5-fp}

\Author{Boris A. KUPERSHMIDT~$^\dag$ and Ognyan S. STOYANOV~$^\ddag$}

\Adress{$\dag$~Department of Mathematics,  University of Tennessee Space Institute,\\
~~Tullahoma, TN 37388,  USA\\
~~E-mail: bkupersh@utsi.edu\\[2mm]
$\ddag$~Department of Mathematics,  Rutgers University,  New Brunswick, NJ 08903, USA\\
~~E-mail: stoyanov@math.rutgers.edu}

\Date{Received February 10, 1999; Accepted March 10, 1999}

\begin{abstract}
\noindent
A Poisson-Lie group acting by the coadjoint action on the dual of its Lie algebra induces
on it a non-trivial class of quadratic Poisson structures extending the linear Poisson bracket
on the coadjoint orbits.
\end{abstract}

\section{Introduction}

If $G$ is a Lie group with Lie algebra ${\mathcal G}$ then the coadjoint action of $G$ on the
dual space ${\mathcal G}^*$ to ${\mathcal G}$ leaves invariant the  linear Poisson bracket on
${\mathcal G}^*$. In other words, the action map $Ad^*:G\times {\mathcal G}^*\to{\mathcal G}^*$
is a Poisson map, with the Poisson structure on $G$ being the trivial one. If $G$ is a Poisson-Lie
group then this action, in general, is no longer Poisson unless the Poisson structure on
${\mathcal G}^*$ could be suitably modif\/ied. In this paper we show that this is indeed possible.
We construct this extension explicitly and give necessary and suf\/f\/icient conditions for its
existence. In particular, any Poisson structure on a f\/inite-dimensional connected
simply-connected Poisson-Lie group $G$ coming from an $r$-matrix that satisf\/ies the Classical
Yang-Baxter Equation induces a Poisson structure on ${\mathcal G}^*$ such that the
coadjoint action becomes Poisson.
The existence of a modif\/ication of the linear Poisson bracket for the case $G=GL(n)$ and
${\mathcal G}^*=gl(n)^*$ was shown earlier in~\cite{Ku1:paper} (see also the related
article~\cite{MaFrei:paper}). In~\cite{Ku1:paper} an equivariant quantization of the coadjoint
action for the case $G=GL(2)$ and ${\mathcal G}^*=gl(2)^*$ was also constructed.

\section{Main Theorem}

We start with recalling some basics. Let $G$ be a f\/inite-dimensional connected
simply-connected Poisson-Lie group. Let the Poisson structure on $G$ be
given by the skew-symmetric rank-2 contravariant tensor ${\pi}^{ij}$. We write the group
law $G\times G {\rightarrow} G$ ($(x,y)\mapsto f(x,y)$) in a neighbourhood of the identity as
\[
f^i(x^1,\ldots,x^n,y^1,\ldots,y^n),\qquad i=1,\ldots,n=\dim G,
\]
where $f^i$ are assumed smooth. Let $\varphi :\ G\to G$ be the map of taking an inverse.
We recall the properties of these maps~\cite{Pontr:paper}:
\be
f^i(0,\ldots,0,y^1,\ldots,y^n)=y^i,\label{eqn1}
\ee
\be
f^i(x^1,\ldots,x^n,0,\ldots,0)=x^i,\label{eqn2}
\ee
\be
f^i(x,\varphi(x))=0=f^i(\varphi(x),x).\label{eqn3}
\ee
From these we deduce that
\be
\frac {\partial f^i}{\partial x^j}(0,0)=\delta^i_j=\frac {\partial f^i}{\partial y^j}(0,0),\label{eqn4}
\ee
and
\be
\delta^i_j+\frac {\partial\varphi^i}{\partial x^j}(0)=0.\label{eqn5}
\ee

For $x\in G$ let $L_y(x)=f(y,x)$ and $R_y(x)=f(x,y)$ be the left and right actions by an
element $y\in G$. The adjoint action $\mbox{Ad}_y(x)=R_{y^{-1}}\circ L_y(x)$ is described in
local  coordinates~by
\be
f^i(f(y,x),\varphi(y)).\label{eqn6}
\ee

Let $I\subset {\mathbb R}$ be an open interval containing the origin and let $x :\ I\to G$ be a
curve $x(t)$ with a direction vector $\xi\in{\mathcal G}=T_{e}G$ at the identity $e=x(0)$.
Dif\/ferentiating~(\ref{eqn5}) in the direction of $\xi$ at $t=0$ (with $x(t)=\exp(t\xi)$) we obtain
a formula for the action of $\mbox{Ad}_y$ on~${\mathcal G}$:
\be
(\mbox{Ad}_y)^i_j\xi^j=\frac {df^i}{dt}\Bigl|_{t=0}=\frac {\partial f^i}{\partial u^k}
(y,\varphi(y))\frac {\partial f^k}{\partial v^j}(y,0)\xi^j\equiv A^i_j(y)\xi^j.\label{eqn7}
\ee
Here $\xi^i$ are the coorsinates of $\xi$ and $u$ and $v$ refer to the f\/irst and second
argument of the functions $f$ respectively. (The Einstein convention of summation over
repeated upper and lower indices is in force throughout this text.) If $\eta\in {\mathcal G}^*$
relative to the form
\be
(\mbox{Ad}_y(\xi),\eta)=(\mbox{Ad}_y)^i_j\xi^j\eta_i\equiv(\xi,\mbox{Ad}_y^*(\eta)),\label{eqn8}
\ee
and the coadjoint action $G\times{\mathcal G}^*\to {\mathcal G}^*$ is def\/ined to be:
\be
(\mbox{Ad}_{y^{-1}}^*)^i_j\eta_i=\frac {\partial f^i}{\partial u^k}(\varphi(y),y)\frac {\partial f^k}
{\partial v^j}(\varphi(y),0)\eta_i=A^i_j\left(y^{-1}\right)\eta_i,\label{eqn9}
\ee
where $\eta_i$ are the coordinates of $\eta$.

\begin{remark} Recall that if $\zeta\in{\mathcal G}$ is a direction vector of a curve $y(s)$ in
$G$, differentiating (\ref{eqn7}) and (\ref{eqn9}) in the direction of $\zeta$ at $s=0$
we obtain formulae for the  maps $\mbox{\rm ad}_{\zeta} :\ {\mathcal G}\to{\mathcal G}$
and $\mbox{\rm ad}_{\zeta}^* :\ {\mathcal G}^*\to{\mathcal G}^*$. For example,
\[\hspace*{-25.55pt}
\ba{l}
\ds (\mbox{\rm ad}_{\zeta}^*)^i_j=\left[-\frac {\partial^2 f^i}{\partial u^k\partial u^s}(0,0)
\frac {\partial f^k}{\partial v^j}(0,0)+\frac {\partial^2 f^i}{\partial u^k\partial v^s}(0,0)
\frac {\partial f^k}{\partial v^j}(0,0)-\frac {\partial^2 f^k}{\partial u^s\partial v^j}(0,0)
\frac {\partial f^i}{\partial u^k}(0,0)\right]\zeta^s
\vspace{3mm}\\
\ds \phantom{(\mbox{\rm ad}_{\zeta}^*)^i_j}=\left[-\frac {\partial^2 f^i}{\partial u^s\partial v^j}(0,0)+
\frac {\partial^2 f^i}{\partial u^j\partial v^s}(0,0)\right]\zeta^s
=-C^i_{sj}\zeta^s,
\ea
\]
and similarly for $\mbox{\rm ad}_{\zeta}$. Here $C^i_{sj}$ are the structure constants of ${\mathcal G}$,
and we have used formulae~(\ref{eqn3}) and (\ref{eqn4}) and
$\frac {\partial^2 f^i}{\partial u^k\partial u^s}(0,0)=0$.
\end{remark}

In addition, the group $G$ is assumed to be equipped with a Poisson-Lie
structure~\cite{Drinfeld:paper1}. That is, a rank-2 contravariant skew-symmetric Poisson
tensor ${\pi}^{kl}$ exists which is compatible with the group multiplication
$G\times G\to G$ (${\pi}$ is a group 1-cocycle). The Poisson bracket between two smooth
functions $f$ and $g$ can be def\/ined in three equivalent ways:
\[
\left\{f,g\right\}={\pi}^{ij}\frac {\partial f}{\partial x^i}\frac {\partial g}{\partial x^j}=
{\varrho}^{ij}{\mathcal D}_{i}f{\mathcal D}_{j}g={\sigma}^{ij}{\mathcal D}^{\prime}_{i}
f{\mathcal D}^{\prime}_{j}g,
\]
where ${\mathcal D}_{i}=dL^{s}_{i}(x)\frac{\partial}{\partial x^s}$ and
${\mathcal D}^{\prime}_{i}=dR^{s}_{i}(x)\frac{\partial}{\partial x^s}$ are left and right
invariant vector f\/ields and $dL(x)$ and $dR(x)$ are the derivatives of the maps
$L_x$ and $R_x$. The tensors ${\pi}$, ${\varrho}$, and ${\sigma}$ are related by
\begin{equation}
{\pi}(x)={\varrho}(x)dL(x)dL(x)={\sigma}(x)dR(x)dR(x).\label{tensor-equivalence}
\end{equation}
The Poisson structure on $G\times G$ is taken to be the product Poisson structure
and thus the multiplication map $m: G\times G\to G$ must be Poisson, that is,
$\left\{f,g\right\}_{G}(xy)=\left\{m^*f,m^*g\right\}_{G\times G}(x,y)$. This implies
 that ${\pi}^{ij}$ must satisfy the 1-cocycle functional equation
\begin{equation}
{\pi}^{ij}(xy)=dR^i_k(y)dR^j_l(y){\pi}^{kl}(x)+dL^i_k(x)dL^j_l(x){\pi}^{kl}(y),
\label{eqn1-cocycle}\end{equation}
which is equivalent to
\begin{equation}
{\sigma}^{ij}(xy)={\sigma}^{ij}(x)+dR^{i}_{k}((yx)^{-1})dR^{j}_{l}((yx)^{-1})A^{k}_{p}(x)
A^{l}_{q}(x)dR^{p}_{s}(yx)dR^{q}_{t}(yx){\sigma}^{st}(y),
\label{eqn2-cocycle}
\end{equation}
where $A(x)=dR(x^{-1})dL(x)=dL(x)dR(x^{-1})$ is the matrix of the adjoint representation.
Let again $y:I\to G$ be a one-parameter curve passing through the identity with direction
vector $\xi$. Dif\/ferentiating (\ref{eqn2-cocycle}) in the direction of $\xi$ at $y=0$
we obtain the system of dif\/ferential equations
\begin{equation}
dL^{s}_{k}(x)\frac{\partial \sigma^{ij}}{\partial x^s}=A^{i}_{s}(x)\alpha^{sp}_{k}A^{j}_{p}(x),
\label{eqn3-cocycle}
\end{equation}
where $\alpha^{ij}_{k}=\frac{\partial \sigma^{ij}}{\partial x^k}(0)=
\frac{\partial \pi^{ij}}{\partial x^k}(0)$ are the components comprising a map
$\alpha :{\mathcal G}\to{\mathcal G}\wedge{\mathcal G}$. The integrability conditions
for (\ref{eqn3-cocycle}) after use of the Maurer-Cartan equations and evaluation at
$x=0$ lead to
\begin{equation}
\alpha^{kl}_sC^s_{ij}=\alpha^{ml}_jC^k_{im}+\alpha^{km}_jC^l_{im}-\alpha^{ml}_i
C^k_{jm}-\alpha^{km}_iC^l_{jm},
\label{eqn4-cocycle}
\end{equation}
and therefore $\alpha :{\mathcal G}\to{\mathcal G}\wedge{\mathcal G}$ being a 1-cocycle
is a necessary and suf\/f\/icient condition for the existence of a solution of (\ref{eqn3-cocycle}).
In particular, as is well known \cite{Drinfeld:paper1}, if $\alpha=\delta r$ is a coboundary,
where $r\in{\mathcal G}\wedge{\mathcal G}$, then
\begin{equation}
\alpha^{ij}_n=C^i_{ns}r^{sj}+C^j_{ns}r^{is}.\label{eqn17}
\end{equation}
In this case (\ref{eqn4-cocycle}) is identically satisf\/ied and the equations (\ref{eqn3-cocycle})
can be trivially integrated yielding the solution
\begin{equation}
\sigma^{ij}(x)=A^{i}_{s}(x)r^{sp}A^{j}_{p}(x)+r^{ij}_{0},
\label{eqn3-coboundary-solution}
\end{equation}
where $r^{ij}_{0}$ is a constant skew-symmetric matrix. (A formula of type (\ref{eqn32})
below, with $dL$ instead of $dR$, is used in the proof of this.)
Substituting (\ref{eqn3-coboundary-solution}) back into the functional equation
(\ref{eqn2-cocycle}), and using the fact that the left and right actions commute, we
deduce that it is a solution of the original 1-cocycle equation if and only if $r^{ij}_{0}=-r^{ij}$.
All of the above is standard in the theory of Poisson-Lie groups. We recalled the relevant
facts which we shall need in the sequel.

Our goal is to make the action $\mbox{Ad}^* :\ G\times{\mathcal G}^*\to{\mathcal G}^*$ Poisson.
Therefore we need to construct (covariant) Poisson tensors $\omega_{ij}$ on ${\mathcal G}^*$
compatible with the coadjoint action. Here again the space $G\times{\mathcal G}^*$ is
equipped with the product Poisson structure and thus the map
$\mbox{Ad}^* :\ G\times{\mathcal G}^*\to{\mathcal G}^*$ is required to be Poisson. This
condition is equivalent to the condition  that locally $\omega_{ij}$ must satisfy the following
system of functional equations
\[
\omega_{ij}(B(y,\eta))=\frac {\partial B_i}{\partial \eta_k}\frac {\partial B_j}
{\partial \eta_l}\omega_{kl}(\eta)+\frac {\partial B_i}{\partial y^k}\frac {\partial B_j}
{\partial y^l}{\pi}^{kl}(y),
\]
where we have introduced $B_i(y,\eta):=A^j_i(y^{-1})\eta_j$, or equivalently
\begin{equation}
\omega_{ij}(A(y^{-1})\eta)=A^k_i(y^{-1})A^l_j(y^{-1})\omega_{kl}(\eta)+
\frac {\partial A^s_i}{\partial y^k}(y^{-1})\frac {\partial A^p_j}{\partial y^l}(y^{-1})
\eta_s\eta_p{\pi}^{kl}(y).\label{eqn10}
\end{equation}
In order to construct solutions we pass to a system of dif\/ferential equations which is
the inf\/initesimal part of (\ref{eqn10}).
Dif\/ferentiation of the above equations in the directions of the coordinate axes yields
\be
\ba{l}
\ds \frac{\partial\omega_{ij}}{\partial\eta_s}\frac{\partial A^l_s}{\partial y^n}\eta_l=
\frac{\partial A^k_i}{\partial y^n}A^l_j\omega_{kl}(\eta)+A^k_i\frac{\partial A^l_j}{\partial y^n}
\omega_{kl}(\eta)+\frac{\partial A^s_i}{\partial y^k}
\frac{\partial A^p_j}{\partial y^l}\frac{\partial {\pi}^{kl}}{\partial y^n}\eta_s\eta_p
\vspace{3mm}\\
\ds \phantom{\frac{\partial\omega_{ij}}{\partial\eta_s}\frac{\partial A^l_s}{\partial y^n}\eta_l=}
+\frac{\partial^2 A^s_i}{\partial y^k\partial y^n}\frac{\partial A^p_j}{\partial y^l}{\pi}^{kl}(y)
\eta_s\eta_p+\frac{\partial A^s_i}{\partial y^k}
\frac{\partial^2 A^p_j}{\partial y^l\partial y^n}{\pi}^{kl}(y)\eta_s\eta_p.
\ea\label{eqn11}
\ee
Evaluating at the identity $y=0$ we obtain
\begin{equation}
C^l_{sn}\eta_l\frac{\partial\omega_{ij}}{\partial\eta_s}=C^k_{in}\omega_{kj}+
C^l_{jn}\omega_{il}-C^s_{ik}C^p_{jl}\alpha^{kl}_n\eta_s\eta_p,\label{eqn12}
\end{equation}
where $\alpha^{kl}_n=\frac{\partial {\pi}^{kl}}{\partial y^n}(0)$ and we have used formulae
$A^i_j(0)=\delta^i_j$, $\frac{\partial A^i_j}{\partial y^n}(0)=C^i_{jn}$, and ${\pi}^{kl}(0)=0$.
We now seek solutions of the above system of dif\/ferential equations.

\begin{remark}
If the Poisson tensor ${\pi}$ on $G$ were zero, then (\ref{eqn12}) reduces to
\begin{equation}
C^l_{sn}\eta_l\frac{\partial\omega_{ij}}{\partial\eta_s}=
C^k_{in}\omega_{kj}+C^l_{jn}\omega_{il}.\label{eqn13}
\end{equation}
It is immediate that the linear bracket $\omega_{ij}(\eta)=C^s_{ij}\eta_s$ satisf\/ies (\ref{eqn13}).
Indeed, after substituting into (\ref{eqn13}) we obtain
\[
\left(C^s_{ln}C^l_{ij}+C^l_{ni}C^s_{lj}+C^l_{jn}C^s_{li}\right)\eta_s=0,
\]
which is identically satisf\/ied. Moreover, any tensor of the form $\omega_{ij}(\eta)=
C^s_{ij}\eta_s\Theta(\eta)$ is a solution of (\ref{eqn13}) as long as $\Theta$ is a solution
of the system
\begin{equation}
C^l_{sn}C^k_{ij}\eta_l\eta_k\frac{\partial\Theta}{\partial\eta_s}=0.\label{eqn14}
\end{equation}
\end{remark}

The general solution of (\ref{eqn12}) is a linear combination of the general solution of the
homogeneous system (\ref{eqn13}) and a particular solution. We look for a particular
solution in the form
\begin{equation}
\omega_{ij}(\eta)=\beta^{qr}_{ij}\eta_q\eta_r,\label{eqn15}
\end{equation}
where $\beta^{qr}_{ij}=-\beta^{rq}_{ji}$ are symmetric in the upper and skew-symmetric
in the lower indices. Substituting into (\ref{eqn12}) we obtain
\[
\left[C^q_{sn}\beta^{sr}_{ij}+C^r_{sn}\beta^{qs}_{ij}-C^s_{in}\beta^{qr}_{sj}-C^s_{jn}\beta^{qr}_{is}-
\frac{1}{2}\left(C^q_{is}C^r_{jl}+C^r_{is}C^q_{jl}\right)\alpha^{sl}_n\right]\eta_q\eta_r=0,
\]
which leads to
\begin{equation}
C^q_{sn}\beta^{sr}_{ij}+C^r_{sn}\beta^{qs}_{ij}-C^s_{in}\beta^{qr}_{sj}-
C^s_{jn}\beta^{qr}_{is}-\frac{1}{2}\left(C^q_{is}C^r_{jl}+C^r_{is}C^q_{jl}\right)\alpha^{sl}_n=0.
\label{eqn16}
\end{equation}

\begin{remark}
The most general tensor $\beta^{qr}_{ij}$ symmetric in the upper and skew-symmetric
in the lower indices that one can construct out of the tensors $\alpha^{ij}_n$ and $C^k_{ij}$ is
\begin{equation}
\beta^{kl}_{ij}=a\left[\alpha^{ks}_iC^l_{sj}-\alpha^{ks}_jC^l_{si}+\alpha^{ls}_i
C^k_{sj}-\alpha^{ls}_jC^k_{si}\right],
\label{eqn_beta}
\end{equation}
where $a$ is a constant scalar. However, (\ref{eqn_beta}) falls short of satisfying (\ref{eqn16})
by the term
\begin{equation}
-aC^m_{ij}\left[\alpha^{qs}_nC^r_{ms}+\alpha^{rs}_nC^q_{ms}\right],
\label{eqn_falls_short}
\end{equation}
where $a=-1/4$. The tensors $\alpha^{ij}_n$ and $C^k_{ij}$ are related by (\ref{eqn4-cocycle}).

In the case when $\alpha$ is a coboundary and given by (\ref{eqn17}) equation (\ref{eqn16})
reads
\begin{equation}
C^q_{sn}\beta^{sr}_{ij}+C^r_{sn}\beta^{qs}_{ij}-C^s_{in}\beta^{qr}_{sj}-C^s_{jn}\beta^{qr}_{is}-
\frac{1}{2}\! \left(C^q_{is}C^r_{jl}+C^r_{is}C^q_{jl}\right)
\!\! \left(C^s_{np}r^{pl}+C^l_{np}r^{sp}\right)\!=0,\label{eqn18}
\end{equation}
and (\ref{eqn_beta}) reduces to
\begin{equation}
\beta^{kl}_{ij}=\frac{1}{2}\left(C^k_{ip}C^l_{js}+C^l_{ip}C^k_{js}\right)r^{sp}+
\frac{1}{4}C^s_{ij}\left[C^l_{sm}r^{km}+C^k_{sm}r^{lm}\right].\label{eqn18a}
\end{equation}
\end{remark}
It turns out that only the f\/irst half in the above formula yields a solution of (\ref{eqn18}).

\begin{proposition} The tensor
\begin{equation}
\beta^{kl}_{ij}=\frac{1}{2}\left(C^k_{ip}C^l_{js}+C^l_{ip}C^k_{js}\right)r^{sp}.\label{eqn19}
\end{equation}
  is a solution of (\ref{eqn18}). Moreover,
$\omega_{ij}(\eta)=\beta^{qr}_{ij}\eta_q\eta_r$ is a Poisson tensor, that is, it satisfies the
Jacobi identities
\begin{equation}
\omega_{ij}\frac {\partial\omega_{kl}}{\partial
\eta_i}+\omega_{ik}\frac {\partial\omega_{lj}}{\partial
\eta_i}+\omega_{il}\frac {\partial\omega_{jk}}{\partial \eta_i}=0,\label{eqn19a}
\end{equation}
if and only if $r$ is a solution of the Classical Yang-Baxter Equation:
\begin{equation}
C^{n}_{sp}r^{sj}r^{pl}+C^{j}_{sp}r^{sl}r^{pn}+C^{l}_{sp}r^{sn}r^{pj}=0,\label{eqn19b}
\end{equation}
provided that a special linear map from ${\mathcal G}^{\otimes 3}$ to
${\mathcal G}^{\otimes 3}\otimes{\mathcal G}^{*\otimes 3}$ has a zero kernel.
\end{proposition}
{\it Proof:} The proof is a straightforward calculation. First, substituting $\beta^{kl}_{ij}$
 from (\ref{eqn19}) into the left hand side of (\ref{eqn18}) and rearanging terms we obtain:
\[
\ba{l}
\left[C^{s}_{ip}C^{q}_{sn}+C^{s}_{ni}C^{q}_{sp}+C^{s}_{pn}C^{q}_{si}\right]C^{r}_{jl}r^{lp}+
\left[C^{s}_{ip}C^{r}_{sn}+C^{s}_{ni}C^{r}_{sp}+C^{s}_{pn}C^{r}_{si}\right]C^{q}_{jl}r^{lp}
\vspace{2mm}\\
\qquad + \left[C^{s}_{jl}C^{r}_{sn}+C^{s}_{nj}C^{r}_{sl}+C^{s}_{ln}C^{r}_{sj}\right]C^{q}_{ip}r^{lp}+
\left[C^{s}_{jl}C^{q}_{sn}+C^{s}_{nj}C^{q}_{sl}+C^{s}_{ln}C^{q}_{sj}\right]C^{r}_{ip}r^{lp},
\ea
\]
which is identically equal to zero. To prove that $\omega_{ij}(\eta)=\beta^{qr}_{ij}\eta_q\eta_r$
is Poisson we note that the Jacobi identities (after symmetrization) are equivalent to the
following identitties for the components of $\beta$:
\begin{equation}
\beta^{qr}_{si}\beta^{sm}_{jk}+\beta^{rm}_{si}\beta^{sq}_{jk}+\beta^{mq}_{si}\beta^{sr}_{jk}+
\mbox{cyclic}(i,j,k)=0.\label{eqn20}
\end{equation}
After a lengthy calculation, with $\beta$ given by (\ref{eqn19}), and using only the fact
that the r-matrix is skew-symmetric, $r^{kl}=-r^{lk}$, and the identities
\begin{equation}
C^{s}_{ip}C^{q}_{sn}+C^{s}_{ni}C^{q}_{sp}+C^{s}_{pn}C^{q}_{si}=0\label{eqn21}
\end{equation}
for the structure constants $C^{k}_{ij}$ of the group, we obtain from (\ref{eqn20}) the equations:
\begin{equation}
\ba{l}
\left[C^{q}_{js}\left(C^{m}_{ku}C^{r}_{iw}+C^{r}_{ku}C^{m}_{iw}\right)+\mbox{cyclic}(q,m,r)\right]
\vspace{2mm}\\
\qquad \times \left[C^{s}_{vt}r^{uv}r^{wt}+C^{w}_{vt}r^{sv}r^{ut}+C^{u}_{vt}r^{wv}r^{st}\right]=0.
\ea
\label{eqn22}
\end{equation}
Let us sketch the major steps of the calculation. Substitution of (\ref{eqn19}) into the left
hand side of (\ref{eqn20}) results in a sum of 36 terms which we group into 9 groups
each consisting of~4 summands:
{\advance\topsep-3pt
\be
C^{q}_{sv}C^{s}_{jt}[C^{r}_{iu}C^{m}_{kw}r^{uv}r^{wt}+C^{m}_{iu}C^{r}_{kw}r^{uv}r^{wt}-
C^{r}_{kw}C^{m}_{iu}r^{wv}r^{ut}-C^{m}_{kw}C^{r}_{iu}r^{wv}r^{ut}]\label{group1}
\ee
\be
\quad +C^{q}_{sv}C^{s}_{kw}[C^{r}_{iu}C^{m}_{jt}r^{uv}r^{wt}+C^{m}_{iu}C^{r}_{jt}r^{uv}r^{wt}-
C^{r}_{jt}C^{m}_{iu}r^{uw}r^{vt}-C^{m}_{jt}C^{r}_{iu}r^{uw}r^{vt}]\label{group2}
\ee
\be
\quad+C^{r}_{sv}C^{s}_{jt}[C^{q}_{iu}C^{m}_{kw}r^{uv}r^{wt}+C^{m}_{iu}C^{q}_{kw}r^{uv}r^{wt}-
C^{q}_{kw}C^{m}_{iu}r^{wv}r^{ut}-C^{m}_{kw}C^{q}_{iu}r^{uv}r^{ut}]\label{group3}
\ee
\be
\quad+C^{r}_{sv}C^{s}_{kw}[C^{q}_{iu}C^{m}_{jt}r^{uv}r^{wt}+C^{m}_{iu}C^{q}_{jt}r^{uv}r^{wt}-
C^{q}_{jt}C^{m}_{iu}r^{uw}r^{vt}-C^{m}_{jt}C^{q}_{iu}r^{uw}r^{vt}]\label{group4}
\ee
\be
\quad+C^{m}_{sv}C^{s}_{jt}[C^{r}_{iu}C^{q}_{kw}r^{uv}r^{wt}+C^{q}_{iu}C^{r}_{kw}r^{uv}r^{wt}-
C^{r}_{kw}C^{q}_{iu}r^{wv}r^{ut}-C^{q}_{kw}C^{r}_{iu}r^{wv}r^{ut}]\label{group5}
\ee
\be
\quad+C^{m}_{sv}C^{s}_{kw}[C^{r}_{iu}C^{q}_{jt}r^{uv}r^{wt}+C^{q}_{iu}C^{r}_{jt}r^{uv}r^{wt}-
C^{r}_{jt}C^{q}_{iu}r^{tv}r^{wu}-C^{q}_{jt}C^{r}_{iu}r^{tv}r^{wu}]\label{group6}
\ee
\be
\quad+C^{q}_{sv}C^{s}_{iw}[C^{r}_{ju}C^{m}_{kt}r^{uv}r^{wt}+C^{m}_{ju}C^{r}_{kt}r^{uv}r^{wt}-
C^{r}_{kt}C^{m}_{ju}r^{tv}r^{wu}-C^{m}_{kt}C^{r}_{ju}r^{tv}r^{wu}]\label{group7}
\ee
\be
\quad+C^{r}_{sv}C^{s}_{iw}[C^{q}_{ju}C^{m}_{kt}r^{uv}r^{wt}+C^{m}_{ju}C^{q}_{kt}r^{uv}r^{wt}-
C^{q}_{kt}C^{m}_{ju}r^{tv}r^{wu}-C^{m}_{kt}C^{q}_{ju}r^{tv}r^{wu}]\label{group8}
\ee
\be
\quad+C^{m}_{sv}C^{s}_{iw}[C^{r}_{ju}C^{q}_{kt}r^{uv}r^{wt}+C^{q}_{ju}C^{r}_{kt}r^{uv}r^{wt}-
C^{q}_{ju}C^{r}_{kt}r^{tv}r^{wu}-C^{q}_{kt}C^{r}_{ju}r^{tv}r^{wu}].\label{group9}
\ee
Each of the above 9 expressions is further transformed to an expression of only 2 summands.
We describe in detail how this is done for the f\/irst of the above groups~(\ref{group1}).
Using the Jacobi identities for the structure constants of ${\mathcal G}$
\begin{equation}
C^{q}_{sv}C^{s}_{jt}+C^{q}_{sj}C^{s}_{tv}+C^{q}_{st}C^{s}_{vj}=0\quad
\Longleftrightarrow\quad C^{q}_{sv}C^{s}_{jt}=-C^{q}_{sj}C^{s}_{tv}-C^{q}_{st}C^{s}_{vj},
\end{equation}
we transform (\ref{group1}) to
\[\hspace*{-14.16pt}
\left(C^{q}_{sj}C^{s}_{tv}+C^{q}_{st}C^{s}_{vj}\right)\left[-C^{r}_{iu}C^{m}_{kw}r^{uv}r^{wt}-C^{m}_{iu}
C^{r}_{kw}r^{uv}r^{wt}+C^{r}_{kw}C^{m}_{iu}r^{wv}r^{ut}+C^{m}_{kw}C^{r}_{iu}r^{wv}r^{ut}\right]
\]
\be
=C^{q}_{sj}C^{m}_{kw}\left[-C^{s}_{tv}C^{r}_{iu}r^{uv}r^{wt}+C^{s}_{tv}C^{r}_{iu}r^{wv}r^{ut}\right]
\label{subgroup1}
\ee
\be
\qquad +C^{q}_{sj}C^{r}_{kw}\left[-C^{s}_{tv}C^{m}_{iu}r^{uv}r^{wt}+C^{s}_{tv}C^{m}_{iu}r^{wv}r^{ut}
\right]\label{subgroup2}
\ee
\be
\qquad +C^{r}_{iu}C^{m}_{kw}\left[-C^{q}_{st}C^{s}_{vj}r^{uv}r^{wt}+C^{q}_{st}C^{s}_{vj}r^{wv}r^{ut}
\right]\label{subgroup3}
\ee
\be
\qquad +C^{r}_{kw}C^{m}_{iu}\left[-C^{q}_{st}C^{s}_{vj}r^{uv}r^{wt}+C^{q}_{st}C^{s}_{vj}r^{wv}r^{ut}
\right].\label{subgroup4}
\ee
Each of the terms (\ref{subgroup1})--(\ref{subgroup4}) is now transformed as follows:
for (\ref{subgroup1}) we have
\[
C^{q}_{sj}C^{m}_{kw}\left[-C^{s}_{tv}C^{r}_{iu}r^{uv}r^{wt}+C^{s}_{tv}C^{r}_{iu}r^{wv}r^{ut}\right]
\]
\[
\qquad =C^{q}_{sj}C^{m}_{kw}\left[-C^{s}_{tv}C^{r}_{iu}r^{uv}r^{wt}+C^{s}_{vt}C^{r}_{iu}r^{uv}r^{wt}
\right]
\]
\[
\qquad =C^{q}_{sj}C^{m}_{kw}\left[C^{s}_{vt}C^{r}_{iu}r^{uv}r^{wt}+C^{s}_{vt}C^{r}_{iu}r^{uv}r^{wt}
\right]
\]
\be
\qquad =2C^{q}_{sj}C^{m}_{kw}C^{s}_{vt}C^{r}_{iu}r^{uv}r^{wt};\label{subgr1}
\ee
for (\ref{subgroup2}) we have
\[
C^{q}_{sj}C^{r}_{kw}\left[-C^{s}_{tv}C^{m}_{iu}r^{uv}r^{wt}+C^{s}_{tv}C^{m}_{iu}r^{wv}r^{ut}\right]
\]
\[
\qquad =C^{q}_{sj}C^{r}_{kw}\left[-C^{s}_{tv}C^{m}_{iu}r^{uv}r^{wt}+C^{s}_{vt}C^{m}_{iu}r^{uv}r^{wt}
\right]
\]
\[
\qquad =C^{q}_{sj}C^{r}_{kw}\left[C^{s}_{vt}C^{m}_{iu}r^{uv}r^{wt}+C^{s}_{vt}C^{m}_{iu}r^{uv}r^{wt}
\right]
\]
\be
\qquad =2C^{q}_{sj}C^{r}_{kw}C^{m}_{iu}C^{s}_{vt}r^{uv}r^{wt};\label{subgr2}
\ee
for (\ref{subgroup3}) we have
\[
C^{r}_{iu}C^{m}_{kw}\left[-C^{q}_{st}C^{s}_{vj}r^{uv}r^{wt}+C^{q}_{st}C^{s}_{vj}r^{wv}r^{ut}\right]
\]
\[
\qquad =C^{r}_{iu}C^{m}_{kw}\left[-C^{q}_{st}C^{s}_{vj}r^{uv}r^{wt}+C^{q}_{sv}C^{s}_{tj}r^{wt}r^{uv}
\right]
\]
\[
\qquad =C^{r}_{iu}C^{m}_{kw}\left[-C^{q}_{st}C^{s}_{vj}-C^{q}_{sv}C^{s}_{jt}\right]r^{uv}r^{wt}
\]
\be
\qquad =-C^{r}_{iu}C^{m}_{kw}C^{q}_{sj}C^{s}_{vt}r^{uv}r^{wt};\label{subgr3}
\ee
and f\/inally for (\ref{subgroup4}) we have
\[
C^{r}_{kw}C^{m}_{iu}\left[-C^{q}_{st}C^{s}_{vj}r^{uv}r^{wt}+C^{q}_{st}C^{s}_{vj}r^{wv}r^{ut}\right]
\]
\[
\qquad =C^{r}_{kw}C^{m}_{iu}\left[-C^{q}_{st}C^{s}_{vj}r^{uv}r^{wt}+C^{q}_{sv}C^{s}_{tj}r^{uv}r^{wt}
\right]
\]
\[
\qquad =C^{r}_{kw}C^{m}_{iu}\left[-C^{q}_{st}C^{s}_{vj}-C^{q}_{sv}C^{s}_{jt}\right]r^{uv}r^{wt}
\]
\be
\qquad =-C^{r}_{kw}C^{m}_{iu}C^{q}_{sj}C^{s}_{vt}r^{uv}r^{wt}.\label{subgr4}
\ee
Adding together (\ref{subgr1})--(\ref{subgr4}) we obtain for (\ref{group1}) the following expression
\[
(\ref{group1})=C^{q}_{sj}C^{s}_{vt}C^{m}_{kw}C^{r}_{iu}r^{uv}r^{wt}+C^{q}_{sj}C^{s}_{vt}
C^{r}_{kw}C^{m}_{iu}r^{uv}r^{wt}.
\]
Performing analogous manipulations as above with each of the terms (\ref{group2})--(\ref{group9})
we deduce that the left hand side of (\ref{eqn20}) is equivalent to
\be
C^{q}_{sj}C^{s}_{vt}C^{m}_{kw}C^{r}_{iu}r^{uv}r^{wt}+C^{q}_{sj}C^{s}_{vt}
C^{r}_{kw}C^{m}_{iu}r^{uv}r^{wt}\label{gro1}
\ee
\be
\qquad +C^{q}_{sk}C^{s}_{vw}C^{r}_{iu}C^{m}_{jt}r^{uv}r^{wt}+C^{q}_{sk}C^{s}_{vw}
C^{m}_{iu}C^{r}_{jt}r^{uv}r^{wt}\label{gro2}
\ee
\be
\qquad +C^{r}_{sj}C^{s}_{vt}C^{m}_{kw}C^{q}_{iu}r^{uv}r^{wt}+C^{r}_{sj}C^{s}_{vt}
C^{m}_{iu}C^{q}_{kw}r^{uv}r^{wt}\label{gro3}
\ee
\be
\qquad +C^{r}_{sk}C^{s}_{vw}C^{q}_{iu}C^{m}_{jt}r^{uv}r^{wt}+C^{r}_{sk}
C^{s}_{vw}C^{m}_{iu}C^{q}_{jt}r^{uv}r^{wt}\label{gro4}
\ee
\be
\qquad +C^{m}_{sj}C^{r}_{iw}C^{s}_{tv}C^{q}_{ku}r^{uv}r^{wt}+C^{m}_{sj}C^{q}_{iw}
C^{s}_{tv}C^{r}_{ku}r^{uv}r^{wt}\label{gro5}
\ee
\be
\qquad +C^{m}_{sk}C^{q}_{jt}C^{s}_{vw}C^{r}_{iu}r^{uv}r^{wt}+C^{m}_{sk}C^{q}_{iu}
C^{s}_{vw}C^{r}_{jt}r^{uv}r^{wt}\label{gro6}
\ee
\be
\qquad +C^{q}_{si}C^{r}_{ju}C^{s}_{vw}C^{m}_{kt}r^{uv}r^{wt}+C^{q}_{si}C^{m}_{ju}
C^{s}_{vw}C^{r}_{kt}r^{uv}r^{wt}\label{gro7}
\ee
\be
\qquad +C^{r}_{si}C^{q}_{ju}C^{s}_{vw}C^{m}_{kt}r^{uv}r^{wt}+C^{r}_{si}C^{m}_{ju}
C^{s}_{vw}C^{q}_{kt}r^{uv}r^{wt}\label{gro8}
\ee
\be
\qquad +C^{m}_{si}C^{r}_{ju}C^{s}_{vw}C^{q}_{kt}r^{uv}r^{wt}+C^{m}_{si}C^{q}_{ju}
C^{s}_{vw}C^{r}_{kt}r^{uv}r^{wt}\label{gro9},
\ee
where each of the pair of terms in (\ref{gro1})--(\ref{gro9}) are obtained from the quadruples
of terms in (\ref{group1})--(\ref{group9}) correspondingly. We now rearange the above
18 terms in the following 6 groups each consisting of 3 terms:
\be
C^{q}_{sj}C^{s}_{vt}C^{m}_{kw}C^{r}_{iu}r^{uv}r^{wt}+C^{m}_{sk}C^{q}_{jt}C^{s}_{vw}
C^{r}_{iu}r^{uv}r^{wt}+C^{r}_{si}C^{q}_{ju}C^{s}_{vw}C^{m}_{kt}r^{uv}r^{wt}\label{grou1}
\ee
\be
\qquad +C^{q}_{sj}C^{s}_{vt}C^{r}_{kw}C^{m}_{iu}r^{uv}r^{wt}+C^{r}_{sk}C^{s}_{vw}
C^{m}_{iu}C^{q}_{jt}r^{uv}r^{wt}+C^{m}_{si}C^{q}_{ju}C^{s}_{vw}C^{r}_{kt}r^{uv}r^{wt}\label{grou2}
\ee
\be
\qquad +C^{r}_{sj}C^{s}_{vt}C^{m}_{kw}C^{q}_{iu}r^{uv}r^{wt}+C^{m}_{sk}C^{q}_{iu}C^{s}_{vw}
C^{r}_{jt}r^{uv}r^{wt}+C^{q}_{si}C^{r}_{ju}C^{s}_{vw}C^{m}_{kt}r^{uv}r^{wt}\label{grou3}
\ee
\be
\qquad +C^{r}_{sj}C^{s}_{vt}C^{m}_{iu}C^{q}_{kw}r^{uv}r^{wt}+C^{q}_{sk}C^{s}_{vw}C^{m}_{iu}
C^{r}_{jt}r^{uv}r^{wt}+C^{m}_{si}C^{r}_{ju}C^{s}_{vw}C^{q}_{kt}r^{uv}r^{wt}\label{grou4}
\ee
\be
\qquad +C^{q}_{sk}C^{s}_{vw}C^{r}_{iu}C^{m}_{jt}r^{uv}r^{wt}+C^{m}_{sj}C^{r}_{iw}C^{s}_{tv}
C^{q}_{ku}r^{uv}r^{wt}+C^{r}_{si}C^{m}_{ju}C^{s}_{vw}C^{q}_{kt}r^{uv}r^{wt}\label{grou5}
\ee
\be
\qquad +C^{r}_{sk}C^{s}_{vw}C^{q}_{iu}C^{m}_{jt}r^{uv}r^{wt}+C^{m}_{sj}C^{q}_{iw}C^{s}_{tv}
C^{r}_{ku}r^{uv}r^{wt}+C^{q}_{si}C^{m}_{ju}C^{s}_{vw}C^{r}_{kt}r^{uv}r^{wt}\label{grou6}.
\ee
The 3 terms in each of the above 6 groups we manipulate further. We show the steps
for~(\ref{grou1}). Thus, for (\ref{grou1}) we have
\[
C^{q}_{sj}C^{s}_{vt}C^{m}_{kw}C^{r}_{iu}r^{uv}r^{wt}+C^{m}_{sk}C^{q}_{jt}C^{s}_{vw}
C^{r}_{iu}r^{uv}r^{wt}+C^{r}_{si}C^{q}_{ju}C^{s}_{vw}C^{m}_{kt}r^{uv}r^{wt}
\]
\[
\qquad =C^{q}_{sj}C^{s}_{vt}C^{m}_{kw}C^{r}_{iu}r^{uv}r^{wt}+C^{m}_{wk}C^{q}_{js}C^{w}_{vt}
C^{r}_{iu}r^{uv}r^{ts}+C^{r}_{ui}C^{q}_{js}C^{u}_{vt}C^{m}_{kw}r^{sv}r^{tw}
\]
\[
\qquad =C^{q}_{sj}C^{s}_{vt}C^{m}_{kw}C^{r}_{iu}r^{uv}r^{wt}+C^{m}_{kw}C^{q}_{sj}C^{w}_{vt}
C^{r}_{iu}r^{uv}r^{ts}+C^{r}_{iu}C^{q}_{sj}C^{u}_{vt}C^{m}_{kw}r^{sv}r^{tw}
\]
\[
\qquad =C^{q}_{sj}C^{m}_{kw}C^{r}_{iu}\left[C^{s}_{vt}r^{uv}r^{wt}+C^{w}_{vt}r^{uv}r^{ts}+
C^{u}_{vt}r^{sv}r^{tw}\right]
\]
\[
\qquad =C^{q}_{sj}C^{m}_{kw}C^{r}_{iu}\left[C^{s}_{vt}r^{uv}r^{wt}+C^{w}_{vt}r^{sv}r^{ut}+
C^{u}_{vt}r^{wv}r^{st}\right]
\]
\[
\qquad =C^{q}_{sj}C^{m}_{ku}C^{r}_{iw}\left[C^{s}_{vt}r^{wv}r^{ut}+C^{u}_{vt}r^{sv}r^{wt}+
C^{w}_{vt}r^{uv}r^{st}\right]
\]
\be
\qquad =C^{q}_{js}C^{m}_{ku}C^{r}_{iw}\left[C^{s}_{vt}r^{wt}r^{uv}+C^{u}_{vt}r^{st}r^{wv}+
C^{w}_{vt}r^{ut}r^{sv}\right].\label{gr1}
\ee
In a completely similar way we obtain for (\ref{grou2})--(\ref{grou6}) the following expressions:
\be
(\ref{grou2})=C^{q}_{js}C^{m}_{iw}C^{r}_{ku}\left[C^{s}_{vt}r^{wt}r^{uv}+C^{u}_{vt}r^{st}r^{wv}+
C^{w}_{vt}r^{ut}r^{sv}\right],\label{gr2}
\ee
\be
(\ref{grou3})=C^{r}_{js}C^{m}_{ku}C^{q}_{iw}\left[C^{s}_{vt}r^{wt}r^{uv}+C^{u}_{vt}r^{st}r^{wv}+
C^{w}_{vt}r^{ut}r^{sv}\right],\label{gr3}
\ee
\be
(\ref{grou4})=C^{r}_{js}C^{m}_{iw}C^{q}_{ku}\left[C^{s}_{vt}r^{wt}r^{uv}+C^{u}_{vt}r^{st}r^{wv}+
C^{w}_{vt}r^{ut}r^{sv}\right],\label{gr4}
\ee
\be
(\ref{grou5})=C^{m}_{js}C^{q}_{ku}C^{r}_{iw}\left[C^{s}_{vt}r^{wt}r^{uv}+C^{u}_{vt}r^{st}r^{wv}+
C^{w}_{vt}r^{ut}r^{sv}\right],\label{gr5}
\ee
\be
(\ref{grou6})=C^{m}_{js}C^{q}_{iw}C^{r}_{ku}\left[C^{s}_{vt}r^{wt}r^{uv}+
C^{u}_{vt}r^{st}r^{wv}+C^{w}_{vt}r^{ut}r^{sv}\right].\label{gr6}
\ee
Finally, adding (\ref{gr1})--(\ref{gr6}) we obtain (\ref{eqn22}).
} 

Thus, we arrive at the following dichotomy. The ``if'' part in the statement of the proposition
follows immediately from (\ref{eqn22}). The ``only if'' part follows from (\ref{eqn22}), provided
the linear map ${\mathcal C} : {\mathcal G}^{\otimes 3}\to{\mathcal G}^{\otimes 3}\otimes
{\mathcal G}^{*\otimes 3}$ with matrix components
\begin{equation}
\ba{l}
{\mathcal C}^{qmr}_{sj,kw,iu}=C^{q}_{js}\left(C^{m}_{ku}C^{r}_{iw}
+C^{r}_{ku}C^{m}_{iw}\right)
\vspace{2mm}\\
\phantom{{\mathcal C}^{qmr}_{sj,kw,iu}=}
+C^{r}_{js}\left(C^{m}_{ku}C^{q}_{iw}+C^{q}_{ku}C^{m}_{iw}\right)+
C^{m}_{js}\left(C^{q}_{ku}C^{r}_{iw}+C^{r}_{ku}C^{q}_{iw}\right),
\ea\label{eqn23}
\end{equation}
has a zero kernel. We do not know how to interpret geometrically this condition on
the group $G$. On the other hand, if $G$ is such that ${\mathcal C}^{qmr}_{sj,kw,iu}=0$, then
any skew-symmetric matrix $r$ induces a Poisson structure on ${\mathcal G}^{*}$, given
by (\ref{eqn19}). This concludes the proof. \hfill \rule{3mm}{3mm}

Two Poisson tensors $\omega^{(1)}_{ij}$ and $\omega^{(2)}_{ij}$ are said to form
a {\it Poisson pair} if their linear combination $a\omega^{(1)}_{ij}+b\omega^{(2)}_{ij}$ is also
a Poisson tensor for arbitrary constants $a$ and $b$. It is easy to see that the Poisson
tensors $\omega^{(1)}_{ij}$ and $\omega^{(2)}_{ij}$ form a Poisson pair if and only if
\begin{equation}
\omega^{(1)}_{ij}\frac {\partial\omega^{(2)}_{kl}}{\partial
\eta_i}+\omega^{(1)}_{ik}\frac {\partial\omega^{(2)}_{lj}}{\partial
\eta_i}+\omega^{(1)}_{il}\frac {\partial\omega^{(2)}_{jk}}{\partial \eta_i}+\omega^{(2)}_{ij}
\frac {\partial\omega^{(1)}_{kl}}{\partial
\eta_i}+\omega^{(2)}_{ik}\frac {\partial\omega^{(1)}_{lj}}{\partial
\eta_i}+\omega^{(2)}_{il}\frac {\partial\omega^{(1)}_{jk}}{\partial \eta_i}=0.\label{eqn24}
\end{equation}
\begin{proposition}
The Poisson tensors
\[
\omega^{(1)}_{ij}=C^s_{ij}\eta_s\Theta(\eta) \qquad and \qquad
\omega^{(2)}_{ij}=1/2\left(C^k_{ip}C^l_{js}+C^l_{ip}C^k_{js}\right)r^{sp}\eta_k\eta_l
\]
form a Poisson pair.
\end{proposition}
{\it Proof:} After substituting $\omega^{(1)}_{ij}$ and $\omega^{(2)}_{ij}$ in (\ref{eqn24})
and collecting terms we obtain
\begin{equation}
\ba{l}
\ds \left(C^{s}_{uj}C^{u}_{kp}+C^{s}_{uk}C^{u}_{pj}+C^{s}_{up}C^{u}_{jk}\right)
C^{v}_{lq}r^{qp}\eta_{s}\eta_{v}
\vspace{3mm}\\
\ds \qquad \qquad +C^{v}_{jq}r^{qp}\eta_{v}\left[C^{u}_{ip}C^{s}_{kl}\eta_{u}\eta_{s}
\frac{\partial \Theta}{\partial \eta_{i}}\right]+\mbox{cyclic}(j,k,l)=0.
\ea \label{eqn25}
\end{equation}
From (\ref{eqn14}) and (\ref{eqn21}) follows that this is an identity. \hfill \rule{3mm}{3mm}

We can thus summarize the above results in the following theorem.
\begin{theorem}
For any finite-dimensional connected simply connected Poisson-Lie group there exists the
family of Poisson structures
\begin{equation}
\omega_{ij}(\eta)=C^s_{ij}\eta_s\Theta(\eta)+C^k_{ip}C^l_{js}r^{sp}\eta_k\eta_l,\label{eqn30}
\end{equation}
on the dual of its Lie algebra such that it makes the coadjoint action Poisson.
Here $\Theta$ is an arbitrary invariant function on ${\mathcal G}^*$.
\end{theorem}
{\it Proof:} What remains to be proved is that the solution of the inf\/initesimal part of~(\ref{eqn10})
obtained above is actually an invariant Poisson bracket under the coadjoint action of the group.
In other words we need to show that the tensor (\ref{eqn30}) satisf\/ies the functional
equation~(\ref{eqn10}). For this we rewrite equation~(\ref{eqn10}) in a new equivalent form.
This is done by using the following properties of the Lie group $G$ and the map $\mbox{Ad}$.
Let $x:I\to G$ be a curve passing through the identity of $G$. We have
\begin{equation}
A^{i}_{j}(f(x(t),y))=A^{i}_{l}(x(t))A^{l}_{j}(y).
\label{eqn31}
\end{equation}
Dif\/ferentiating at $t=0$ we obtain
\begin{equation}
\frac{\partial A^{i}_{j}}{\partial y^k}(y)dR^{k}_{s}(y)=C^{i}_{ls}A^{l}_{j}(y)\quad
\Longrightarrow\quad \frac{\partial A^{i}_{j}}{\partial y^k}(y)=C^{i}_{ls}A^{l}_{j}(y)dR^{s}_{k}(y^{-1}).
\label{eqn32}
\end{equation}
From the identity $A^{i}_{s}(y)A^{s}_{j}(y^{-1})=\delta^{i}_{j}$ we have also (after dif\/ferentiating
in the coordinate directions)
\begin{equation}
\frac{\partial A^{i}_{j}}{\partial y^k}(y^{-1})=-A^{i}_{s}(y^{-1})\frac{\partial A^{s}_{p}}{\partial y^k}(y)
A^{p}_{j}(y^{-1}).
\label{eqn33}
\end{equation}
We note also the invariance of the constant tensor $C^{i}_{jk}$:
\begin{equation}
C^{i}_{jk}=A^{i}_{s}(y)C^{s}_{pq}A^{p}_{j}(y^{-1})A^{q}_{k}(y^{-1}).
\label{eqn34}
\end{equation}
Now we substitute (\ref{eqn33}) into the original functional equation (\ref{eqn10}) and
after an easy calculation using (\ref{eqn32}) and (\ref{eqn34}) we transform it to the equivalent
equation
\begin{equation}
\ba{l}
\omega_{ij}(A(y^{-1})\eta)=A^k_i(y^{-1})A^l_j(y^{-1})\omega_{kl}(\eta)
\vspace{2mm}\\
\qquad \qquad + C^{m}_{iq}C^{n}_{jt}A^s_m(y^{-1})A^p_n(y^{-1})\eta_s\eta_pdR^{q}_{k}
(y^{-1})dR^{t}_{l}(y^{-1}){\pi}^{kl}(y).
\ea \label{eqn35}
\end{equation}
Using the relation (\ref{tensor-equivalence}) between the tensors $\pi$ and $\sigma$ we f\/inally
obtain
\begin{equation}
\omega_{ij}(A(y^{-1})\eta)=A^k_i(y^{-1})A^l_j(y^{-1})\omega_{kl}(\eta)+C^{m}_{iq}C^{n}_{jt}
A^s_m(y^{-1})A^p_n(y^{-1})\eta_s\eta_p{\sigma}^{qt}(y).
\label{eqn36}\end{equation}
With $\sigma$ given by formula (\ref{eqn3-coboundary-solution}) it is now a straightforward
calculation to verify that the tensor $\omega$ as given by formula (\ref{eqn30}) satisf\/ies the
functional equation (\ref{eqn36}). Indeed, the left hand side of (\ref{eqn36}) after substitution
of (\ref{eqn30}) reads
\begin{equation}
C^s_{ij}A^{p}_{s}(y^{-1})\eta_p\Theta(A(y^{-1})\eta)-C^s_{iu}C^p_{jv}r^{uv}A^{m}_{s}
(y^{-1})A^{n}_{p}(y^{-1})\eta_m\eta_n.
\label{eqn37}
\end{equation}
The right hand side of (\ref{eqn36}) yields
\be
\ba{l}
A^{k}_{i}(y^{-1})A^{l}_{j}(y^{-1})C^s_{kl}\eta_s\Theta(\eta)-A^{k}_{i}(y^{-1})A^{l}_{j}(y^{-1})
C^s_{ku}C^p_{lv}r^{uv}\eta_s\eta_p
\vspace{2mm}\\
+C^m_{iq}C^n_{jt}A^{s}_{m}(y^{-1})A^{p}_{n}(y^{-1})\eta_s\eta_pA^{q}_{u}(y)r^{uv}A^{t}_{v}(y)
-C^m_{iq}C^n_{jt}A^{s}_{m}(y^{-1})A^{p}_{n}(y^{-1})\eta_s\eta_pr^{qt}.
\ea \label{eqn38}
\ee
Using
\begin{equation}
C^m_{iq}A^{q}_{u}(y)A^{s}_{m}(y^{-1})\stackrel{(\ref{eqn34})}{=}C^s_{qu}A^{q}_{i}(y^{-1}),\qquad
C^n_{jt}A^{t}_{v}(y)A^{p}_{n}(y^{-1})\stackrel{(\ref{eqn34})}{=}C^p_{mv}A^{m}_{j}(y^{-1}),
\label{eqn39}\end{equation}
we transform the third term of (\ref{eqn38}) and after comparison of terms in the left and right
hand sides we conclude that the functional equation (\ref{eqn36}) is identically satisf\/ied.
Thus, every solution $r$ of the Classical Yang-Baxter Equation induces the Poisson
structure (\ref{eqn30}) on ${\mathcal G}^*$ making it a homogeneous Poisson space under
the coadjoint action of $G$. This concludes the proof.\hfill \rule{3mm}{3mm}

\section{Discussion}

Here we show that the newly obtained Poisson bracket on ${\mathcal G}^*$ when specialized
to the case of $G=GL(n)$ and ${\mathcal G}^*=gl(n)^*$ recovers the one obtained
in \cite{Ku1:paper}. Let $x^{\alpha}_i$ be the components of a matrix representation
$x : GL(n)\to \mbox{Mat}(n)$ of $GL(n)$. The multiplication map $(x,y)\mapsto f(x,y)$ is given by
\begin{equation}
f_i^{\alpha}(x,y)=x^{\alpha}_sy^s_{i}.
\label{eqn40}
\end{equation}
Then it is easy to compute the structure constants
\begin{equation}
C^{\alpha j k}_{i \beta \gamma}=\left[\frac {\partial^2 f^{\alpha}_i}
{\partial u^{\beta}_j\partial v^{\gamma}_k}(0,0)-
\frac {\partial^2 f^{\alpha}_i}{\partial u^{\gamma}_k\partial v^{\beta}_j}(0,0)\right]
=\delta^{\alpha}_{\beta}\delta^{j}_{\gamma}\delta^{k}_{i}
-\delta^{\alpha}_{\gamma}\delta^{k}_{\beta}\delta^{j}_{i}.
\label{eqn41}
\end{equation}
With these structure constants and the r-matrix components
$r^{\alpha\beta}_{ij}=-r^{\beta\alpha}_{ji}$ we have
\[\hspace*{-5.51pt}
\left\{\eta^{j}_{\beta},\eta^{k}_{\gamma}\right\}=C^{\alpha j k}_{i \beta \gamma}\eta^{i}_{\alpha}
\Theta(\eta)+C^{\alpha k m}_{i \gamma\omega}C^{\epsilon j n}_{l \beta\lambda}
\eta^{i}_{\alpha}\eta^{l}_{\epsilon}r^{\lambda\omega}_{nm}
\]
\[\hspace*{-5.51pt}
\phantom{\left\{\eta^{j}_{\beta},\eta^{k}_{\gamma}\right\}}
=(\delta^{\alpha}_{\beta}\delta^{j}_{\gamma}\delta^{k}_{i}-\delta^{\alpha}_{\gamma}\delta^{k}_{\beta}
\delta^{j}_{i})\eta^{i}_{\alpha}\Theta(\eta)+(\delta^{\alpha}_{\gamma}
\delta^{k}_{\omega}\delta^{m}_{i}-\delta^{\alpha}_{\omega}\delta^{m}_{\gamma}
\delta^{k}_{i})(\delta^{\epsilon}_{\beta}\delta^{j}_{\lambda}\delta^{n}_{l}
-\delta^{\epsilon}_{\lambda}\delta^{n}_{\beta}\delta^{j}_{l})
\eta^{i}_{\alpha}\eta^{l}_{\epsilon}r^{\lambda\omega}_{nm}
\]
\[\hspace*{-5.51pt}
\phantom{\left\{\eta^{j}_{\beta},\eta^{k}_{\gamma}\right\}}
=(\delta^{j}_{\gamma}\eta^{k}_{\beta}-\delta^{k}_{\beta}\eta^{j}_{\gamma})\Theta(\eta)+
(\delta^{k}_{\omega}\eta^{m}_{\gamma}-\delta^{m}_{\gamma}\eta^{k}_{\omega})
(\delta^{j}_{\lambda}\eta^{n}_{\beta}-\delta^{n}_{\beta}\eta^{j}_{\lambda})
r^{\lambda\omega}_{nm}
\]
\be\hspace*{-5.51pt}
\phantom{\left\{\eta^{j}_{\beta},\eta^{k}_{\gamma}\right\}}
=(\delta^{j}_{\gamma}\eta^{k}_{\beta}-\delta^{k}_{\beta}\eta^{j}_{\gamma})\Theta(\eta)+
r^{\lambda\omega}_{\beta\gamma}\eta^{k}_{\omega}\eta^{j}_{\lambda}+
r^{jk}_{nm}\eta^{n}_{\beta}\eta^{m}_{\gamma}-r^{\lambda k}_{\beta m}
\eta^{m}_{\gamma}\eta^{j}_{\lambda}-r^{j\omega}_{n\gamma}\eta^{k}_{\omega}\eta^{n}_{\beta},
\label{eqn42}
\ee
which is in agreement with the formula obtained in \cite{Ku1:paper}.

Since there is a canonical isomorphism $T^*G\simeq G\times{\mathcal G}^*$ and the Poisson
 tensor $\omega^{(2)}_{ij}$ on ${\mathcal G}^*$ forms a Poisson pair with the linear tensor
$\omega^{(1)}_{ij}$ on ${\mathcal G}^*$, the corresponding tensors on $T^*G$  will also form
a Poisson pair under this isomorphism. It will be interesting to construct quantizations of the
cotangent space $T^*G\simeq G\times{\mathcal G}^*$ and the group $G$ and lift the coadjoint
Poisson action to an equivariant quantum action between non-commutative spaces in the
quantum case. We hope to address this problem in a future publication.

\label{kupershmidt_5-lp}

\end{document}